\documentstyle[12pt,amscd]{amsart}

\newcommand{\Cal}{\cal}
\newcommand{\nc}{\Bbb{C}}
\newcommand{\g}{\frak{g}}

\newtheorem{th}{Theorem}
\newtheorem{prop}[th]{Proposition}

\theoremstyle{remark}

\newtheorem{ack}{Acknowledgment}  

\begin{document}
\title{Cohomology and deformation of Leibniz pairs}
\author{M. Flato}
\address{Physique Math\'ematique, Universit\'e de Bourgogne,
B.P. 138, F-21004 Dijon Cedex, France}
\email{flato@@satie.u-bourgogne.fr}

\author{M. Gerstenhaber}
\address{
Department of Mathematics, University of Pennsylvania,
Philadelphia, PA 19104-6395, USA}
\email{mgersten@@mail.sas.upenn.edu
 or murray@@math.upenn.edu (TeX files)}
\thanks{The second author wishes to thank the U.S. National
Security Agency for partial support during the preparation of
this paper}

\author{A. A. Voronov}
\address{Department of Mathematics, University of Pennsylvania,
Philadelphia, PA 19104-6395, USA}
\email{voronov@@math.upenn.edu}
\thanks{
The third author wishes to thank the National
Science Foundation for partial support during
 the preparation of
this paper}

\date{December 26, 1994}

\begin{abstract}
Cohomology and deformation theories are developed for Poisson
algebras starting with the more general concept of a Leibniz
pair, namely of an associative algebra $A$ together with a Lie
algebra $L$ mapped into the derivations of $A$. A bicomplex
(with both Hochschild and Chevalley-Eilenberg cohomologies) is essential.
\end{abstract}

\maketitle

The importance of Poisson algebras in classical mechanics makes
it useful to have a deformation theory for them. To construct
this we are led to define a more general concept: A {\bf Leibniz pair}
$(A,L)$ consists of an associative algebra $A$  and a Lie algebra
$L$ over some common coefficient ring $k,$ connected by a Lie
algebra morphism $\mu:L \to \operatorname{Der}A,$ the Lie algebra
of $k$-linear derivations of  $A$ into itself.  It can happen that $L$ is
identical with $A$ as a $k$-module; but note that in the graded case
discussed below the gradings may differ.  When this is so, denoting
the associative product of $a,b \in A$ by $ab$ and their Lie
product by $[a,b]$, if one has further that $\mu (a)(bc) = [a,bc] =
[a,b]c + b[a,c]$, then
$A$ will be called a {\bf non-commutative Poisson
algebra} (``ncPa'').  A Poisson algebra in the usual sense is one
where the associative multiplication on $A$ is commutative (and
$k = \Bbb R$ or $\Bbb C$). A  Leibniz pair with $A = L$ need not
be a Poisson algebra since one may have, for example, an abelian
Lie multiplication together with a non-trivial structure morphism
$\mu.$

This notion of a non-commutative Poisson algebra is a particular case
of that suggested by P. Xu \cite{ping}. In his definition, he
considers an associative algebra $A$ together with a class $\Pi \in
H^2(A, A)$ such that $[\Pi, \Pi] = 0$ in $H^3 (A, A)$ in the sense of
the Gerstenhaber bracket on Hochschild cohomology of $A$. Xu's
definition is especially suitable for the geometric situation, such as
$A = \nc[x_1,
\dots, x_n]$ or $A = C^\infty (M)$, when $M$ is a manifold and one
takes multilinear differential operators as cochains. Then $H^2(A, A)$
will be naturally isomorphic to the space of bivector fields and the
condition $[\Pi, \Pi] = d B$, where, due Hodge theory for Hochschild
cohomology, one can assume $B$ is symmetric; this will
yield the Jacobi identity for $\Pi$.  Thus, $\Pi$ will give rise to a
Poisson bracket on $A$. For a general commutative or noncommutative
algebra $A$, $H^2(A,A)$ would not be exhausted by skew biderivations
of $A$, and the relevance of a class $\Pi$ to any kind of Poisson
bracket on $A$ would be lost. Our definition translates into this
language as follows. A non-commutative Poisson algebra is nothing but
an associative algebra $A$ and a skew biderivation $\Pi \in C^2(A, A)$
satisfying the Jacobi identity. It will not even automatically be a
2-cocycle, unless $A$ is commutative.

A Leibniz pair $(A,L)$ in which $A$ is commutative, $L$ is an $A$-module, and
the morphism $\mu:L \to \operatorname{Der}A$ is an $A$-module morphism will be
called a {\bf Rinehart pair}, \cite{rinehart, palais, YKSM},
provided abbreviating
\[
\mu (x)a = \mu (x,a)  = [x,a]: L \otimes A \to A,
\]
with $x \in L, a \in A, $
the following natural condition is also satisfied
\[
[x, ay] = [x,a] y + a[x,y].
\]
Here the commutativity of $A$ as an
associative algebra is essential, for it is only in this case that the
derivations of $A$ form a module over $A.$ Any commutative Poisson algebra $A$
gives rise to a Rinehart pair $(A,A)$ but one must be careful in notation:
writing $\operatorname{ad}b$ for the mapping sending $c$ to $[b,c]$ one must
distinguish between $a \operatorname{ad}b$ and $\operatorname{ad} ab.$

Leibniz pairs $(A,L)$ are common but the concept seems not to have been
previously considered in this generality. The simplest example is that where
$A$ is an arbitrary associative algebra and $L$ is just $\operatorname{Der} A$
or $A$ itself, considered as the Lie algebra of inner derivations of $A$. This
includes the case where $A$ is either an ordinary polynomial ring
$k[x_1,\dots,x_n]$ or free non-commutative polynomial ring $k\langle
x_1,\dots,x_n\rangle$ in $n$ variables and $L$ the ``linear first-order
differential operators'', i.e., the free $k$-module spanned by $\partial
/\partial x_1,\dots,\partial/\partial x_n.$ (In the second case $A$ is
identical with the tensor algebra over $k$ on the module spanned over $k$ by
the variables).

A most important example of a commutative Leibniz pair is
$(C^{\infty}(\Cal M), \operatorname{Vect} \Cal M)$, where $A$ is the
algebra of smooth functions on a manifold $\Cal M$ and $L$ the Lie algebra of
smooth vector fields; this is a Rinehart pair.  Of special interest is the case
when ${\Cal M}$ is symplectic or Poisson (for a recent monograph on
Poisson manifolds see \cite{vaisman}),
$A = C^\infty({\Cal M})$ and $L = A$ endowed
with the Poisson bracket (central extension by the constants of the Lie
algebra of Hamiltonian vector fields); there the theory of star
products \cite {BFFLS}, associative deformations of $A$ giving rise to
Lie algebra deformations of $L$, is a prototype of the theory we
are developing here; in particular the interplay between the
(infinite) Hochschild cohomology of $A$ and the (finite) Chevalley-Eilenberg
cohomology of $L$ was essential there, and is in fact one of the main
motivations for introducing Leibniz pairs and Poisson algebras.

Leibniz pairs also arise whenever a Lie group $G$ operates
as a group of automorphisms on an associative
algebra $A;$ the pair is then $(A,\frak g)$, where $\frak g =
\operatorname{Lie} G,$ the Lie algebra of $G.$  (One may view the Lie group as
acting on a ``non-commutative function space'' $A.$)
This situation occurs frequently in mathematical physics when $A$
is a $C^*$ algebra of observables (crossed products of a $C^*$ algebra
and a group of automorphisms), including the so-called theory of KMS states
(generalized equilibrium states) where ${\frak g}$ is one-dimensional,
generating a semi-group. A special case arises in
geometry when a Lie group $G$ acts on a manifold $\Cal M$.
The pair is then $(C^{\infty}(\Cal M), \frak g)$.

Analogous to the ``canonical'' Leibniz
pair $(A, \operatorname{Der}A)$ associated to an associative
algebra $A$, if one has a Lie algebra $L$ then one can form the
Leibniz pair $(\Cal UL,L),$ where $\Cal UL$ is the universal
enveloping algebra of $L.$

Another interesting example arises in the context of quantum groups at roots of
unity. Suppose $\Cal U_q (\g)$ is a quantum group. Specializing $q$ to a root
$\epsilon$ of 1, one usually gets a big center $Z$. The center is then a
Poisson algebra with respect to the bracket
\[
\{x,y\}  =  \lim_{q \to \epsilon} \frac{[\tilde x, \tilde y] }{q- \epsilon},
\]
where $\tilde x$ and $\tilde y \in \Cal U_q (\g)$ denote liftings of
elements $x$ and $y \in Z$. This bracket can be extended to an action
of $Z$ on $\Cal U_\epsilon (\g)$ by derivations, and thus $(\Cal
U_\epsilon (\g), Z)$ is another example of a Leibniz pair.

In this paper we develop cohomology and deformation theories for Leibniz pairs,
Rinehart pairs, ncPa's, and Poisson algebras. For the latter it generally will
not matter if a Poisson algebra is commutative or not,  so ``Poisson algebra''
will henceforth include the non-commutative case unless it is said otherwise.
The theory will be written in an algebraic context, but its extension to a
topological context in the infinite case (at least for Fr\'echet nuclear
algebras, as in \cite {BFGP}) is straightforward.

As observed a Leibniz pair may be viewed as an infinitesimal version
of an algebra with a group of operators. The deformation theory of
such objects is clear if the group is kept fixed \cite{gersts}. It
seems less natural, however, to consider simultaneous deformation of
the algebra (in some sense, a local object) together with that of the
group (a global one). Passage to the Leibniz pair (which consists of
two local objects) may be a remedy. In suitable cases the passage from
an algebra with operators to a Leibniz pair $(A,L)$ can be reversed,
for example, when $A$ and $L$ are both finite-dimensional (real or
complex), the adjoint group of $L$ then operates on $A$.

\section{Note on the graded case and operads}

While we do not explicitly consider the graded case here, our methods remain
applicable. An important open question in both the ungraded and graded cases is
the relationship between our cohomology theory and that constructed by other
methods.  Poisson algebras are special cases of operad algebras, cf.\
\cite{gk}, for which there is already a general cohomology theory. Further
Getzler and Jones
\cite{gj} have studied a cohomology theory for Gerstenhaber algebras, one of
the two natural graded generalizations of Poisson algebras, which we
will describe in this Section. In the more general context of algebras
with many operations grouped into dot-like and bracket-like ones,
satisfying some distributivity laws, Fox and Markl \cite{fmarkl} have
studied cohomology using the operad cohomology viewpoint as well. Fox
and Markl show that for such algebras, the complex computing
cohomology is the total complex of a bicomplex, cf.\ Section
\ref{double} below.
 Similar considerations are found in the cyclic cohomology bicomplex
 \cite {Connes}, where also one direction is Hochschild cohomology and which
 classifies the so-called \cite{CFS} closed star products (for which a trace
can be
 defined in the star product algebra).

We consider first {\bf graded Leibniz pairs} $(A,L)$.  Here both $A$ and $L$
are graded, degrees of elements being indicated by superscripts, and $L$ is a
graded Lie algebra, cf.\ \cite{gerst}. That is,
\begin{gather*}
[x^m,y^n] = -(-1)^{mn}[y^n,x^m] \quad \text{and} \\
(-1)^{mp}[x^m,[y^n,z^p]]+(-1)^{nm}[y^n,[z^p,x^m]]
+(-1)^{pn}[z^p,[x^m,y^n]] = 0.
\end{gather*}
The first equation asserts that the multiplication in $L$ is
``graded skew'', the second is the graded Jacobi identity. Such
$L$ are often called ``Lie superalgebras.'' A graded
derivation of degree $p$ of $A$ is a linear map $\Cal D:A \to A$
with $\Cal D(A^m) \subset A^{m+p},$ (where $A^m$ denotes the
homogeneous part of $A$ of degree $m$) and
$$\Cal D(a^mb^n) = (\Cal D a^m)b^n + (-1)^{pm}a^m(\Cal D b^n).$$
These form a graded Lie algebra, and the structure morphism
$\mu:L \to \operatorname{Der} A$ is now required to be a morphism
of graded Lie algebras, hence to carry $z^p$ to a homogeneous
derivation of $A$ of degree $p.$

The associative algebra $A$ is called graded commutative or
``supercommutative'' if $a^mb^n = (-1)^{mn}b^na^m,$  as, for example, in a
Grassmann algebra. There are now two significantly different generalizations of
an ordinary Poisson algebra, both of which have their ``right'' and ``left''
versions.

In a right ``Poisson'' superalgebra, $A$ is simultaneously a graded
commutative algebra and a graded Lie algebra {\it with the same grading}, and
{\it right} Lie multiplication by $c^p$ ({\it i.e.,} the operation carrying $a$
to $[a,c^p]$) is a graded derivation of degree $p.$  That is, we have
$[a^mb^n,c^p] = [a^m,c^p]b^n + (-1)^{mp}a^m[b^n,c^p].$ In the left version, it
is {\it left} Lie multiplication by $c^p$ which is a graded derivation of
degree $p$, i.e., we have $[c^p, a^mb^n] = [c^p,a^m]b^n +
(-1)^{mp}a^m[c^p,b^n].$ There is no significant difference between these, since
one can pass from one to the other by replacing the associative multiplication
in $A$ by its ``opposite'', i.e., by reversing its order, but keeping the Lie
operation unchanged.  A right  Poisson superalgebra then becomes a left one,
and vice versa.

In a right {\bf Gerstenhaber algebra}, $A$ and $L$ are
identical as {\it ungraded} $k$-modules, but as graded modules
the Lie degree is reduced by 1, so one has
\begin{gather*}
 [a^m,b^n] = -(-1)^{(m-1)(n-1)}[b^n,a^m] \\
\begin{split}
(-1)^{(m-1)(p-1)}[a^m,[b^n,c^p]] +
(-1)^{(n-1)(m-1)}[b^n,[c^p,a^m]] \\
+
(-1)^{(p-1)(n-1)}[c^p,[a^m,b^n]] = 0
\end{split}
\\
[a^mb^n, c^p] =[a^m,c^p]b^n + (-1)^{m(p-1)} a^m[b^n,c^p].
\end{gather*}
The last equation asserts that right Lie multiplication by $c^p$ acts as a
graded derivation of degree $p-1$ of the graded associative algebra $A.$ In a
left Gerstenhaber algebra, it is left Lie multiplication which acts as a graded
derivation.  As before, if the associative multiplication is replaced by its
opposite and the Lie multiplication preserved, then ``right'' and ``left'' are
interchanged.

Gerstenhaber algebras (and in fact the concept of a graded Lie algebra)  first
arose in algebraic deformation theory as the Hochschild cohomology
$H^\bullet(A,A)$
of an associative algebra $A$ with coefficients in itself, \cite{gerst}.
[Under the newer convention which introduces a sign change when an operator
passes an argument, $H^\bullet(A,A)$ becomes a left Gerstenhaber algebra.
Specifically, if $f$ and $g$ are, respectively, $m$-- and $n$--cochains of $A$
with coefficients in itself and $a, b$ respectively $m$-- and $n$--tuples of
elements of $A$, then Hochschild originally defined $(f\cup g)(a,b)$ to be
$f(a)g(b)$ while the newer convention sets $(f\cup g) (a,b) =
(-1)^{mn}f(a)g(b)$ because $g,$ which is viewed as having degree $n$, has
passed $a,$ which has degree $m.$]

Notice that the degree 0 part of a Poisson superalgebra is an ordinary Poisson
algebra, but ``degree zero part'' has less meaning for a Gerstenhaber algebra
since what is of even degree for the associative structure is of odd degree for
the graded Lie structure, and conversely. Nevertheless it is true that
the elements of degree 0 form an associative algebra and those of degree
 1 a Lie algebra.
(Warning: Some authors  unfortunately
use ``Poisson superalgebra '' for both cases.)

\section{A missing equivalence}

Even in the ungraded case, the operad approach to the cohomology of
Poisson algebras does not lend itself easily to exhibiting the basic
relation between cohomology and deformations, and it is not
immediately applicable to Leibniz pairs. The latter problem could be
overcome if there were an equivalence between the category of Leibniz
pairs and a subcategory of the category of associative algebras or
even if Leibniz pairs were equivalent to algebras over a certain
operad. While Leibniz pairs are definitely not themselves algebras
over an operad, it seems likely that the data carried by a Poisson
algebra is equivalent to a number of operations satisfying certain
identities on a single vector space. But this has so far not been
shown. Here we construct a cohomology theory for Leibniz pairs and
algebras in an elementary way, as the cohomology of a certain double
complex, and show that (as in the classical algebraic deformation
theory \cite{g:1964} and that of bialgebras
\cite{gs:1990:PNAS,gs:1992:CM,BFGP}) the first cohomology group contains the
infinitesimal automorphisms, the second cohomology group contains the
infinitesimal deformations and the third group the obstructions to these.
Since the concept of a module over a Leibniz pair is meaningful, as we show in
the next section, other elementary (i.e., non-operadic) definitions of
cohomology are also possible (e.g., as a Yoneda type Ext, or as a derived
functor), but we do not know if all these coincide.  While injective and
projective modules may be defined in the usual way, we do not know if there
exist ``enough'' injective or projective modules over a Leibniz pair $(A,L),$
that is, whether every module is a submodule of an injective or quotient of a
projective. (The existence of enough injectives would imply that the standard
construction of the cohomology using injective resolutions is possible.)

We do not know either if there is a Gerstenhaber algebra structure on the
cohomology of a Leibniz pair with coefficients in itself or even simply a
graded Lie structure. In the special case where the Leibniz pair structure is
induced from a Lie group operation there is some evidence that the Gerstenhaber
algebra structure should exist. We will return to this at the end of the paper.

\section{Modules}

In a Leibniz pair $(A,L)$ we will generally denote the elements of $A$ by $a,
b, c, \dots$ and those of $L$ by $x, y, z, \dots.$ Suppose that $M$ is an
$A$-bimodule and $P$ an $L$--module, and denote by $A+M$ and $L+P$,
respectively, the associative and Lie semidirect products.  Then $(M,P)$ will
be called an $(A,L)$-- {\bf module} if there is a Lie algebra morphism $L +P
\to \operatorname{Der}(A+M)$ extending the structure morphism $\mu:L \to
\operatorname{Der}A$ such that $[L,M], [P,A] \subset M$ and $[P,M] = 0.$ These
are grading conditions: view $L$ and $A$ as of degree zero and $P$ and $M$ as
of degree one; then $[P,M] = 0$  because there is nothing of degree two.  Note
that if $\pi \in P$ then the hypotheses imply that the map of $A$ into itself
sending $a$ to $[\pi,a]$ is a derivation of $A$ into $M$, so we have a map $P
\to \operatorname{Der}(A,M).$

Every Leibniz pair $(A,L)$ is a module over itself, but it is
helpful to distinguish between $A$ the algebra and $A$ the
$A$--bimodule. In this case we may denote the module by $tA$, where
$t$ is some variable whose square is zero.  Then $A+tA$ is
identical with the $k[t]/t^2$--module $A[t]/t^2$ which, when
viewed as a $k$--algebra, is just the semidirect product.  We
will then similarly write $tL.$ If $A$ is a Poisson algebra then
a module $M$ over $A$ is really a module $(M,M)$ over the pair
$(A,A),$ where the two Lie products of  both $M$ and $A$ are
identical.

The concept of a morphism of Leibniz pairs is clear, whence so is
that of a morphism of modules over a Leibniz pair.  One has
kernels and cokernels, hence an abelian category, and also direct
sums in the algebraic sense, i.e, categorical direct products.
The category of Leibniz pairs contains the categories of
associative and Lie algebras as full subcategories.  (In the
first case, take $L = 0,$ in the second $A = k$.)

\section{The double complex for a Leibniz pair}

Unless stated otherwise, all tensor products will be over the
coefficient ring $k.$  (For much of what follows this need not be
of characteristic $0$, but to fix the ideas one may assume that
it is $\Bbb R.$)  Since we do not consider the graded case, we
will henceforth denote $\otimes^pA$ simply by $A^p.$ Suppose that
$(A,L)$ is a Leibniz pair and that $(M,P)$ is a module over it. Let
$C^p(A,M) = \operatorname{Hom}_k(A^p,M),$ the $p$th Hochschild
cochain group of $A$ with coefficients in $M.$  Note that
$C^0(A,M)$ is just $M.$
Now $C^p = C^p(A,M)$ is again an $L$--module in the usual way,
where if $f \in C^p$ then one sets
$$[x,f](a_1,\dots,a_p) = [x,f(a_1,\dots,a_p)]
     - \sum_{i=1}^pf(a_1,\dots,[x,a_i],\dots,a_p).$$
It is immediate that $\operatorname{Der}(A,M)$ is then an
$L$--submodule of $C^1(A,M),$ and the mapping $P \to
\operatorname{Der}(A,M)$ noted above is an $L$--module morphism.

The Hochschild coboundary $\delta_H:C^p \to C^{p+1}$ is given by
\begin{gather*}
(\delta_H f)(a_1,\dots,a_{p+1}) = a_1f(a_2,\dots,a_{p+1})
     +\\
\sum_{i=1}^p(-1)^pf(a_1,\dots,a_ia_{i+1},\dots,a_{p+1})
     +(-1)^{p+1}f(a_1,\dots,a_p)a_{p+1}.
\end{gather*}
It is an essentially trivial matter to verify the following
fundamental result.

\begin{prop} Let $(M,P)$ be a module over a Leibniz pair
$(A,L).$ Then $\delta_H:C^p(A,M) \linebreak[0] \to C^{p+1}(A,M)$ is a morphism
of Lie modules over $L. \qed$
\end{prop}

\begin{sloppypar}
Let $(M,P)$ be a module over the Leibniz pair $(A,L).$ The double
complex $C^{\bullet,\bullet}(A,L;M,P)$ of $(A,L)$ with
coefficients in $(M,P)$ is defined as follows.  For all $p > 0,
q \ge 0,$ set $C^{p,q} = \operatorname{Hom}_k(\wedge^qL,
C^p(A,M)) \cong \operatorname{Hom}_k(A^p\otimes\wedge^qL, M),$
and set
$C^{0,q} = \operatorname{Hom}_k(\wedge^qL, P).$ (Note that $A^0 =
\wedge^0L = k.$)  For $p > 0,$ the ``vertical'' coboundary
$C^{p.q} \to C^{p+1,q}$ is just the Hochschild coboundary
$\delta_H$ (the index $h$ stands for Hochschild, not for horizontal!).
For $p=0$ note that we are given a Lie module
$\operatorname{Hom}_kA = C^1(A,M).$ This induces, for every $q$,
the desired map $C^{0,q} \to C^{1,q}$ which we denote $\delta_v.$
The composition of two vertical coboundary operators  $C^{p,q}
\to C^{p+2,q}$ is zero for all $p$ and $q$  since for $p>0$ both
are Hochschild coboundaries, while for $p=0$ it follows because
$\operatorname{Der}(A,M)$ is precisely the kernel of
$\delta_H:C^1(A,M) \to C^2(A,M).$
\end{sloppypar}

\[
\begin{CD}
\dots && \dots && \dots \\
@A\delta_HAA       @A\delta_HAA  @A\delta_HAA \\
\operatorname{Hom}(A^2, M) @>\delta_{CE}>>
\operatorname{Hom}(A^2\otimes L, M) @>\delta_{CE}>>
\operatorname{Hom}(A^2\otimes\wedge^2L, M) @>\delta_{CE}>> \dots\\
@A\delta_HAA       @A\delta_HAA  @A\delta_HAA \\
\operatorname{Hom}(A, M) @>\delta_{CE}>>
\operatorname{Hom}(A\otimes L, M) @>\delta_{CE}>>
\operatorname{Hom}(A\otimes\wedge^2L, M) @>\delta_{CE}>> \dots\\
@A\delta_vAA       @A\delta_vAA  @A\delta_vAA  \\
\operatorname{Hom}(k, P) @>\delta_{CE}>>
\operatorname{Hom}(L, P) @>\delta_{CE}>>
\operatorname{Hom}(\wedge^2 L, P) @>\delta_{CE}>> \dots
\end{CD}
\]

In the horizontal direction we have for all $p$ and $q$ the
Chevalley-Eilenberg coboundary $\delta_{CE}:C^{p,q} \to
C^{p,q+1}.$ Since all the vertical coboundaries are Lie module
morphisms, it follows that they commute with all the horizontal
ones. Denoting any vertical coboundary for the moment by
$\delta_H$ we can therefore define a total cohomology complex
with $C_{tot}^n = \bigoplus_{p+q=n}C^{p,q}$ and
$\delta_{tot}:C_{tot}^n \to C_{tot}^{n+1}$ whose restriction to
$C^{p,q}$ by definition is $\delta_H + (-1)^p\delta_{CE}.$
Its cohomology we define to be the {\bf cohomology of the
Leibniz pair} $(A,L)$ with coefficients in the module $(M,P),$
denoted $H^{\bullet}(A,L;M,P).$ When we consider the
cohomology of $(A,L)$ with coefficients in itself as a module, we
will write simply $C^{p,q}(A,L)$ and $H^\bullet(A,L).$

Denote the invariants of an $L$--module  $N$, {\it i.e.}, the set
of all $n \in N$ with $[x,n] = 0$ for all $x\in L$ by
$N^{L}.$ In the double complex $C^{\bullet,\bullet}(A,L;M,P)$ the
kernel of $\delta_{CE}:C^{p,0} \to
C^{p,1}$ is just $C^p(A,M)^{L}$ for all $p>0,$ and for $p =
0$ it is $(P)^L.$ (When cohomology is taken with coefficients in
$(A,L)$ itself this will be the center $Z(L)$ of $L.$)  The
``augmenting column" of the double complex is therefore the
complex
$$C^{\bullet}:P^L \to C^1(A,M)^{L} \to C^2(A,M)^{L} \to \dots.$$
The first arrow is induced by the morphism $P \to
\operatorname{Der}(A,M)$ (which, when cohomology is taken with
coefficients in $(A,L)$ is the structural morphism $\mu: L \to
\operatorname{Der}A$) and all the others by the Hochschild
coboundary. Note that $H^i(A,L;M,P)$ is in particular an $L$--
module and is finite-dimensional if all the algebras and modules
are.

While the cohomology groups of the double complex are not easy to
compute, the Whitehead lemmas and the fact that the cohomology of
$L$ with coefficients in itself vanishes give the following

\begin{th}
Suppose that $(M,P)$ is a module over a
Leibniz pair $(A, L)$ with all algebras and modules
finite-dimensional over the real or complex numbers and
$L$ semisimple.  Then $H^i(A, L;M, P)$ is naturally isomorphic to
$H^i(C^{\bullet})$ for $i = 1, 2, 3.  \qed$
\end{th}

If the associative algebra $A$ of a Leibniz pair is commutative,
and if the coefficient ring $k$ has characteristic 0, that is,
if there is a unital morphism $\Bbb Q \to k$ then
each of the $C^p(A,M)$ decomposes into its Hodge parts
\cite{gs:1987} and one may check that this
decomposition is preserved by the Chevalley-Eilenberg coboundary.  The
entire double complex therefore has a Hodge decomposition, the first
component of which is obtained by replacing each $C^p(A,M)$ by its
submodule of Harrison cochains (cochains vanishing on shuffles).  (The
augmenting column similarly decomposes, giving the augmenting columns
of the Hodge components.)

\section{The double complex for a Poisson algebra}
\label{double}

When $A$ is a (not necessarily  commutative) Poisson algebra,
the Lie module operation of $A$ on itself can not be deformed
independently of the Lie multiplication on $A$, so the foregoing
double complex must be modified by the elision of one row.  In
this section we will assume characteristic zero, but for the
deformation theory it is really only necessary to assume the
invertibility of 2 and 3 in $k.$

Let $M$ be a Poisson $A$-module, i.e., simultaneously an
associative $A$--bimodule and a Lie $A$--module where for each
fixed $m \in M$ the mapping sending $a \in A$ to $[a,m] \in M$ is
a derivation of $A$ into $M.$  Note that we do not have to assume
separately that there is given a map $M \to C^1(A,M);$ for this
we simply take the usual Hochschild coboundary in which
$(\delta_Hm)(a) = am - ma.$ (This can vanish identically, in
which case $M$ is called a ``symmetric'' module.  A particular
case is that in which $A$ is a commutative Poisson algebra and $M$
is $A$ itself.)  We will define a double complex $\widetilde
C^{\bullet,\bullet}(A;M)$ with $\widetilde C^{p,q} = C^{p,q} =
\operatorname{Hom}_k((A^p)\otimes(\wedge^qA),M)$ for $p \ge 2$,
but with $\widetilde C^{1,q} =
\operatorname{Hom}_k(\bigwedge^{q+1} A,M)$ and $\widetilde C^{0,q}$
undefined.  The lowest vertical coboundary, from
$\widetilde C^{1,q}$ to $\widetilde C^{2,q},$ will be denoted
$\delta_P,$ (the ``Poisson'' coboundary).
\[
\begin{CD}
\dots && \dots && \dots && \dots \\
@AAA  @A\delta_HAA       @A\delta_HAA  @A\delta_HAA \\
0 @>>> \operatorname{Hom}(A^2, M) @>\delta_{CE}>>
\operatorname{Hom}(A^2\otimes A, M) @>\delta_{CE}>>
\operatorname{Hom}(A^2\otimes\wedge^2 A, M) \to\\
@AAA  @A\delta_PAA       @A\delta_PAA  @A\delta_PAA \\
M @>\delta_{CE}>> \operatorname{Hom}(A, M) @>\delta_{CE}>>
\operatorname{Hom}(\wedge^2 A, M) @>\delta_{CE}>>
\operatorname{Hom}(\wedge^3 A, M) \to
\end{CD}
\]

Suppose now that $V$ is simply a $k$--module and denote
$\otimes^pV$ by $V^p.$  We view $\bigwedge^pV$ as the submodule of $V^p$
consisting of all skew elements and denote by
$\epsilon_p:V^p \to \wedge^pV$ the projection sending
$v_1\otimes\cdots\otimes v_n$ to $v_1\wedge\cdots\wedge v_n.$
If $a \in V^p, b\in V^q$ then
$\epsilon_{p+q}a\otimes b =
\epsilon_{p+q}(\epsilon_pa\otimes\epsilon_qb).$ We will henceforth
write $\epsilon$ without the subscript. Let $\epsilon^*$ be the
dual map $\operatorname{Hom}_k(\wedge^qA,M) \to
\operatorname{Hom}_k(A\otimes(\wedge^{q-1}A),M).$  If $f$ is in
the latter module, we will write $f(a_1|a_2\wedge\cdots\wedge a_q)$
for its value at $a_1\otimes(a_2 \wedge \dots \wedge a_q).$  Now let $\delta_P:
(\bigwedge^qA,M) \to (A^2\otimes \bigwedge A^{q-1})$ be the
composite of $\epsilon^*$ with $\delta_H: \operatorname{Hom}_k
(A\otimes\bigwedge^{q-1}A,M) \to
\operatorname{Hom}_k(A^2\otimes\bigwedge^{q-1} A, M),$ and view it
as a ``vertical'' coboundary.  (Note that this is vertical by
virtue of our renumbering, but that it is really our original row
1 which is missing.)
\begin{prop}
In the original double complex, if $f \in C^{0,q}$ then $\delta_{CE}\epsilon^*
f+ \epsilon^*\delta_{CE}f \linebreak[0]
= \delta_vf$.
\end{prop}
\begin{pf} From the definitions we have
\begin{gather*}
(\delta_{CE}f)(a_1\wedge\cdots\wedge a_{q+1}) =
\sum_{i=1}^{q+1}(-1)^i[a_i,f(a_1\wedge\cdots\wedge\hat a_i\wedge\cdots)] + \\
\sum_{j>i}(-1)^{i+j}f([a_i,a_j] \wedge\cdots\hat a_i\cdots\hat a_j \cdots),
\end{gather*}
while
\begin{gather*}
(\delta_{CE}\epsilon^*f)(a_1|a_2\wedge\cdots\wedge a_{q+1})=
\sum_{i=2}^{q+1}(-1)^i[a_i,(\epsilon^*f)(a_1|a_2\wedge\cdots\wedge \hat a_i
\wedge \cdots a_{q+1})] +\\
\sum_{j>i\ge2}(-1)^{i+j}(\epsilon^*f)(a_1|[a_i,a_j]\wedge\cdots
\hat a_i\cdots\hat a_j\cdots).
\end{gather*}
It follows immediately that
\begin{gather*}
(\epsilon^*\delta_{CE}f)(a_1|a_2\wedge\cdots\wedge a_{q+1}) =\\
-(\delta_{CE}\epsilon^*f)((a_1|a_2\wedge\cdots\wedge a_{q+1}) +\\
[a_1,f(a_2\wedge\cdots\ a_{q+1})] + \sum_{j=2}^{q+1}(-1)^{1+j}
f([a_i,a_j]|a_2\wedge\cdots\hat a_j\wedge\cdots a_{q+1}),
\end{gather*}
which implies the formula.
\end{pf}

It follows, in particular, that in a commutative Poisson algebra
 if $f \in C^{0,q}$ then
$\delta_{CE}\epsilon^* f+ \epsilon^*\delta_{CE}f = 0$.
Note, incidentally, the similarity of the result with the classical formula
in differential geometry $di+id = \Cal L.$

\begin{th}
 If $f \in C^{0,q}$ then
$$\delta_{CE}\delta_Pf + \delta_P\delta_{CE}f = 0.$$
\end{th}
\begin{pf} Apply $\delta_H$ to both sides of the foregoing formula and note
that $\delta_H\delta_v = 0,$ while $\delta_H$ and $\delta_{CE}$ commute.
\end{pf}

It follows that we indeed have a double complex. (Since the horizontal and
vertical coboundaries in our original complex commute, one should really
replace
the horiziontal coboundary on the $q$--th column with $(-1)^q\delta_{CE}$ so
that they now anti-commute. The total coboundary is then just the sum of
horizontal and vertical coboundaries and $\delta_P$ is now the appropriate
first vertical coboundary in the modified double complex.)
There is no evident relation between the cohomology of the double complex for
a Leibniz pair and that for a Poisson algebra. In particular, neither is
included into nor a quotient of the other. Generally, when the bottom rows
of a double complex are altered, the total cohomology groups may change
in unknown ways.

\section{Deformation of Leibniz pairs and Poisson algebras}

The deformation theory of a Leibniz pair now parallels that of an
associative algebra (\cite{g:1964}), the controlling
cohomology being that of the pair with coefficients in itself.

When a Leibniz pair $(A,L)$ is given we will, as before,
generally denote the elements of $A$ by $a,b,\dots$ and of $L$ by
$x,y,\dots$.  Let $\alpha: A \times A \to A$ denote the
associative multiplication map, $\lambda: L \times L \to L$ be
the Lie multiplication, and write $\mu(x)(a) = \mu(x,a)$, viewing
$\mu$ as a map $L \times A \to A.$ The Leibniz pair may then be
denoted $(\alpha, \mu, \lambda).$ By a {\bf deformation} of this
Leibniz pair we will mean a Poisson pair of the form $(\alpha_t,
\mu_t, \lambda_t)$ where $\alpha_t$ is a deformation of $\alpha$
({\it cf}.\ \cite{g:1964}), {\it i.e.,} an associative
$k[[t]]$--bilinear multiplication $A[[t]] \times A[[t]] \to
A[[t]]$ expressible as a power series  $\alpha_t(a,b) = ab +
t\alpha_1(a,b) + t^2\alpha_2(a,b) + \dots $ where each $\alpha_i$
is a bilinear map $A  \times A \to A$ extended to be a
$k[[t]]$--bilinear map $A[[t]] \times A[[t]] \to A[[t]],$ and
similarly for $\lambda_t, \mu_t.$  The triple $(\alpha_1, \mu_1,
\lambda_1)$ of first order terms is the {\bf infinitesimal} of
the deformation. Note that as $\lambda_t$ must be skew, so are
all the $\lambda_i$ and in particular $\lambda_1.$  It is clear
that $\alpha_1$ must be an infinitesimal deformation of $A$ as an
associative algebra, {\it i.e.,} a Hochschild 2-cocycle of $A$
with coefficients in $A$ itself, and that $\lambda_1$ is an
infinitesimal deformation of $L,$ and so a 2-cocycle in the
Chevalley-Eilenberg theory of $L$ with coefficients in $L.$ In
addition, there are compatibility conditions betweeen $\alpha_1$
and $\mu_1$ and between $\mu_1$ and $\lambda_1$ which one obtains
from the linear terms in the following equations, respectively:
\begin{gather*}
\mu_t(x, \alpha_t(a,b)) = \alpha_t(a, \mu_t(x,b) ) +
\alpha_t(\mu_t(a,a),b), \\
\mu_t(\lambda_t(x,y),a) = \lambda_t(x,\mu_t(y,a)) - \mu_t(y,
\mu_t(x,a))
\end{gather*}
These, when written out together with the foregoing, state precisely that
$(\alpha_1,\mu_1,\lambda_1)$ constitute a
2-cocycle of the total cohomology arising from the double complex
$C^{\bullet,\bullet}(A,L;A,L).$ (The idea of the double complex we have
constructed  is borrowed, in spirit, from the construction of the cohomology of
a bialgebra, {\it cf}.\ \cite{gs:1990:PNAS}, or for a better exposition
\cite{gs:1992:CM}.)

\begin{sloppypar}
There is a natural concept of {\bf equivalence of deformations}:
we say that $(\alpha_t, \mu_t, \lambda_t)$  is equivalent to
$(\alpha_t', \mu_t', \lambda_t')$ where $\alpha_t' = \alpha
+t\alpha_1' +t^2\alpha_2' + \dots,$ etc., if there exist
$k[[t]]$--linear maps
$\Phi_t = \operatorname(id)_A+t\phi_1 + t^2\phi_2 + \dots: A[[t]]
\to A[[t]]$
(with each $\phi_i$ a linear map $A \to A$ extended to be
$k[[t]]$--linear) and
$\Psi_t = \operatorname{id}_L +t\psi_1 +t^2\psi_2 + \dots: L[[t]]
\to L[[t]]$ such that
\end{sloppypar}
\begin{gather*}
\alpha_t'(a,b) = \Phi_t^{-1}\alpha_t(\Phi_ta,\Phi_tb), \\
\mu_t'(x,a) = \Phi_t^{-1}\mu_t(\Psi_tx, \Phi_ta), \\
 \lambda_t'(x,y) = \Psi_t^{-1}\lambda_t(\Psi_tx,\Psi_ty).
\end{gather*}
As expected, equivalent deformations have
cohomologous infinitesimals, and  any cocycle in the
cohomology class of the infinitesimal of a deformation can be
taken as the infinitesimal of an equivalent deformation.  Up to
equivalence, therefore, only the cohomology class of an
infinitesimal is of consequence.   Obstructions
arise (at least in principle) exactly as in earlier deformation
theories: if $(\alpha_t, \mu_t, \lambda_t)$ are given through
terms in $t^{n-1}$ and satisfy the definition of a Leibniz pair
modulo $t^n$, then there is generally an obstruction in the third
cohomology group to extending through terms in $t^n$ in such a
way that the conditions for a Leibniz pair are satisfied modulo
$t^{n+1}.$

In the subcategory of Poisson algebras by definition we have $A = L$ and $\mu =
\lambda$, and in the foregoing must take $\Psi_t = \Phi_t.$ We then must use
the double complex $\widetilde C^{\bullet,\bullet}(A;A)$ and have analogous
assertions to the foregoing.  Restriction to the subcategory of commutative
Poisson algebras is then essentially that of the passage from Hochschild to
Harrison cohomology.

In addition to deformations, one-parameter families of automorphisms
(restricted to be the identity at $t = 0$) were considered in
\cite{g:1964}. There an infinitesimal automorphism of an algebra $A$
was simply a derivation, and in characteristic zero there was no
formal obstruction to finding a family with this infinitesimal since
one could exponentiate the derivation. The matter is less simple for a
Leibniz pair $(A,L)$ since an infinitesimal automorphism is now a pair
$(\phi,\psi)$ constituting a one-cocycle of our double complex and it
appears that even in characteristic zero there is an obstruction to
constructing the higher order terms.  For a Poisson algebra $A$, an
infinitesimal automorphism is any $\phi:A \to A$ which is
simultaneously a derivation of both the associative and Lie
structures, and these can be formally exponentiated. Amongst them are
all ad $a$ for $a \in A.$

Recall that when an ordinary associative algebra $A$ with
multiplication $\alpha$ is deformed to $A_t$, the new multiplication
being $\alpha_t$, then it is not always possible to ``lift'' or extend
a cocycle $F$ of $A$ (with coefficients in itself) to one of $A_t,$
cf.\ \cite{g:1974}. The primary obstruction is the class of the graded
Lie product $[F, \alpha_1].$ But even if every derivation could be
lifted, there might be obstructions in the sense of the present
paper. For note that when $(A, \operatorname{Der}A)$ is a Leibniz pair,
a deformation requires that the Lie algebra structure of
$\operatorname{Der}A$ actually be deformed. By contrast, if $A$ alone
is deformed and every element of $\operatorname{Der}A$ is liftable, it
is still conceivable that the Lie algebra which they generate is now
of greater dimension than originally. Also, recall that when a
non-nilpotent Lie group $G$ operates on an associative algebra $A$,
then the various universal deformation formulas based on the group
(cf. \cite{coll-gs:19xx}) generally produce deformations of $A$ which
no longer permit the operation of the full group $G$; symmetry is
generally broken.  Here, however, we are in particular investigating
ways in which $A$ deforms while preserving the full operation of a Lie
group (although the operation itself may simultaneously be deformed).

When $A$ is finite-dimensional, one gets some information about the deformation
problem from a theorem of Kubo, \cite{kubo}. First, $A$ is said to be $L$-{\bf
simple}, if it contains no proper ideal which is simultaneously an
$L$-submodule. Kubo's theorem asserts that such $A$ must be a total matric
algebra,
that is, the algebra of all $n \times n$-matrices for some $n$, or
multiplication in $A$ must be identically zero and $A$ must be nothing more
than an irreducible $L$-module, made into an algebra with the trivial
multiplication. When $L$ is simple, there can be no deformations of this
structure, but for general $L$ the question becomes one of the deformation
theory of modules over a Lie algebra.

\section{Afterword: Application of the Cohomology Comparison Theorem}

For a presheaf of algebras over a small category, one can define in a natural
way a cohomology theory which, remarkably, is identical with the cohomology of
a
single ring built from the presheaf; this is part of the  ``Cohomology
Comparison Theorem'' (CCT) of \cite{gersts}. (One consequence is
that for every simplicial complex there is an algebra whose cohomology with
coefficients in itself is naturally isomorphic to the simplicial cohomology of
the complex.) Since the cohomology has been reduced to that of a single ring,
it follows from the original results of \cite{gerst} that there indeed
exists a Gerstenhaber algebra structure on the cohomology of the original
presheaf.  However, while for a single algebra $A$ there is already a graded
Lie structure on the Hochschild cochains of $A$ with coefficients in itself,
this is not the case for a presheaf of algebras, where the graded Lie structure
exists only on the cohomology level.  It is therefore far from transparent.
Nevertheless, it suggests that when we have a Leibniz pair $(A,L)$ induced
from the action of a Lie group on $A$ then we should apply the CCT to the
group $G$ (considered as a small category) and the presheaf which to the unique
object of this category assigns the algebra $A.$ The essential thing we do not
know is what relation the cohomology so obtained has to that which we have just
defined for the pair $(A,L).$ If these turn out to be the same, then we shall
have, amongst other things, that there is a Gerstenhaber algebra structure on
the cohomology of $(A,L).$ This approach may be applicable to the case where
the given  Leibniz pair consists of the smooth functions and smooth vector
fields on a compact manifold. If it is, then one would expect that the
Gerstenhaber algebra structure is given by explicit (although difficult)
formulas whose nature really should not depend on the fact that we started with
a compact manifold. It would seem, therefore, that there is at least some hope
of exhibiting a Gerstenhaber algebra structure on the cohomology of a Leibniz
pair or Poisson algebra.

\begin{ack}
We are grateful to Y. Kosmann-Schwarzbach, F. Kubo, M. Markl,
J. Stasheff, D. Sternheimer, A. Weinstein, and P. Xu for helpful discussions.
\end{ack}

\makeatletter \renewcommand{\@biblabel}[1]{\hfill#1.}\makeatother
\renewcommand{\bysame}{\leavevmode\hbox to3em{\hrulefill}\,}

\end{document}